\documentclass{elsarticle}
\usepackage{lineno,hyperref}
\modulolinenumbers[5]
\journal{The Art of Discrete and Applied Mathematics}
\usepackage[english]{babel}
\usepackage{amssymb,amsmath,amsfonts}

\title{On divisible design Cayley graphs}
\tnotetext[grant]{
The work is supported by Mathematical Center in Akademgorodok. Both authors are also partially supported by RFBR according to the research project 17-51-560008.
}
\author[01]{Vladislav~V.~Kabanov}
\ead{vvk@imm.uran.ru}
\author[01,02]{Leonid~Shalaginov}
\ead{44sh@mail.ru}

\address[01]{Krasovskii Institute of Mathematics and Mechanics, S. Kovalevskaja st. 16, Yekaterinburg, 620990, Russia}
\address[02]{Chelyabinsk State University, Brat'ev Kashirinyh st. 129, Chelyabinsk,  454021, Russia}

\usepackage{natbib}
\usepackage{graphicx}
\begin{document}
\begin{frontmatter}

\begin{abstract}
Let $v,k,b,a$ be integers such that  $v > k \geq b \geq a \geq 0$. 
A {\em  Deza graph} with parameters $(v,k,b,a)$ is a $k$-regular graph on $v$ vertices  in which the number of common neighbors  of any two distinct vertices takes two values $a$ or $b$ ($a\leq b$).
A $k$-regular graph on $v$ vertices is a {\em divisible design graph} with parameters $(v,k,\lambda_1 ,\lambda_2 ,m,n)$ when its vertex set can be partitioned into $m$ classes of size $n$, such that any two distinct vertices from the same class have  $\lambda_1$ common neighbors, and any two vertices from different classes have  $\lambda_2$ common neighbors. It is clear, that divisible design graphs are  Deza graphs. 

First of all, we proved that divisible design Cayley graphs arise only by means of divisible difference sets  relative to some subgroup. Secondly, we constructed a special set in affine group over finite field and proved that this set is a divisible difference set and thus give us a divisible design graph.

\end{abstract}

\begin{keyword}
Deza graph; divisible design graph; divisible design; divisible different set; Cayley graph; affine group over a finite field.
\vspace{\baselineskip}
\MSC[2010] 	05C75\sep 05B30\sep 05E30
\end{keyword}
\end{frontmatter}

\section{Introduction}

Let $v,k,b,a$ be integers such that  $v > k \geq b \geq a \geq 0$. 
A {\em  Deza graph} with parameters $(v,k,b,a)$ is a $k$-regular graph on $v$ vertices  in which the number of common neighbors  of any two distinct vertices takes two values $a$ or $b$ ($a\leq b$). Deza graphs appeared as a generalization of strongly regular graphs  in   \cite{EFHHH}. This is a wide class of graphs which includes not only strongly regular graphs but also divisible design graphs \cite{HKM, CH}, regular $(0,\lambda)$-graphs  \cite{AB, BO} and others.

A $k$-regular graph on $v$ vertices is a {\em divisible design graph} with parameters $(v,k,\lambda_1 ,\lambda_2 ,m,n)$ when its vertex set can be partitioned into $m$ classes of size $n$, such that any two distinct vertices from the same class have  $\lambda_1$ common neighbors, and any two vertices from different classes have  $\lambda_2$ common neighbors. For a divisible design graph the partition into  classes is called  a {\em canonical partition}.

Let $G$ be a finite group  with the identity element $e$. If $S$ is a subset of $G$ which is closed under inversion and does not contain $e$, then {\em Cayley graph} $\mathbb{C}ay(G,S)$ is a graph with the vertex set $G$ and two vertices $x$, $y$ are adjacent if and only if $x y^{-1} \in S$. 

In this paper, some divisible design graphs are constructed that arise from finite groups in the form of Cayley graphs. The following theorem is the basis for the construction.  
\medskip

{\bf Theorem 1.} {\em  Let $\mathbb{C}ay(G,S)$ be a Deza graph with parameters $(v,k,b,a)$. Let also $A$, $B$ and $\{e\}$ be a partition of $G$ and $S S^{-1}$ be a multiset such that
$S S^{-1} = a A + b B + k\{e\}$. If either $A\cup \{e\}$ or $B\cup \{e\}$ is a subgroup of $G$,  then $\mathbb{C}ay(G,S)$ is a divisible design graph and the right cosets of this subgroup give a canonical partition of the graph. Conversely, if 
$\mathbb{C}ay(G,S)$ is a divisible design graph, then the class of its canonical partition which 
contains the identity of $G$ is a subgroup of $G$ and classes of the canonical partition of divisible design graph coincide with the cosets of this subgroup.}
\medskip

{\bf Proof.} Let $N=A\cup \{e\}$ be a subgroup of $G$ and let $G = N a_1\cup \dots \cup N a_{r}$ be a partition of $G$ by the right cosets of $N$. 
For every $g, h\in N$ and $1\leq i\leq r$, $1\leq j\leq r$, the set of neighbors of $ha_i$ is $Sha_i$ and the set of neighbors of $g a_j$ is $S g a_j$. Thus, the set of common neighbors for $ha_i$ and $g a_j$ coincides with $S h a_i\cap S g a_j$. The number $|S h a_i\cap S g a_j|$ equals $|\{(s, s^{*})| sha_i = s^{*}ga_j\}|$, where $s,s^{*}\in S$. Therefore, $h = s^{-1}s^{*}ga_ja_i^{-1}$. If $i=j$, then $h = s^{-1}s^{*}g$ and $s^{-1}s^{*}\in N$. So there are exactly $a$ such pairs if $s^{-1}s^{*}\in A$ and $k$ such pairs if $s^{-1}s^{*}= e$. If $i\neq j$, then $ga_ja_i^{-1}\notin N$ and hence $s^{-1}s^{*}\notin N$. So $s^{-1}s^{*}\in B$ and there are exactly $b$ such pairs.
The case of $N=B\cup \{e\}$ is viewed in a similar way.

Conversely, let $\mathbb{C}ay(G,S)$ be a divisible design graph and $N$ be a class of the canonical partition of divisible design graph which contains the identity of $G$. It's enough to prove that for any $h,g\in N$ we have
$h g^{-1}\in N$. Since $h$ and $g$ belong to the same class of the canonical partition of  
$\mathbb{C}ay(G,S)$, then $|S h\cap S g|=\lambda_1$.   The number of pairs $(s,s^{*})$ such that $sh=s^{*}g$ is equal to  $\lambda_1$. Thus, $hg^{-1}=s^{-1}s^{*}$ is repeated  $\lambda_1$ times in $S S^{-1}$.
\medskip

Theorem 1 shows that Cayley divisible design graphs arise only by means of divisible difference sets  relative to some subgroup.

Let $G$ be a finite group of order $mn$ and $N$ be a subgroup of $G$ of order $n$. Then  a subset $S$ of $G$ is called a {\em divisible difference set with exceptional subgroup} $N$ if there are constants  $\lambda_1$ and  $\lambda_2$ such that every non-identity element of $N$ can be expressed as a right quotient of elements in $S$ in exactly  $\lambda_1$ ways and
every element in $G \setminus N$ can be expressed as a right quotient of elements in $S$ in exactly  $\lambda_2$ ways.

In other words, if $k=|S|$, then  
$S S^{-1}$ is the following multiset: $$\lambda_1 (N - \{e\}) +  \lambda_2 (G - N) + k\{e\}.$$

\section{Construction of divisible design Cayley graphs}
\medskip

Let $\mathbb F$ be a finite field with $q^r$ elements, where $q$ is a prime power and $r >1$. 

Consider the group $\mathbb{G} $ of all $2\times 2$ matrices 
$\left(\begin{array}{cc}1 & 0\\ \alpha &\beta\\\end{array}\right),$
 where $\alpha\in \mathbb{F}$ and $\beta\in \mathbb{F}\setminus \{0\}$. It's clear that $\mathbb{G}$  is a semi-direct product of two subgroups $N$ and $K$, where $N$ consists of all matrices $\left(\begin{array}{cc}1 & 0\\ \alpha & 1 \\\end{array}\right)$, and $\alpha \in \mathbb{F},$  $K$ consists of all matrices 
$\left(\begin{array}{cc}1 & 0\\ 0 & \beta \\\end{array}\right)$, and $ \beta\in 
\mathbb{F}\setminus \{0\}$. 

 define  a bijection $\psi^{+}$ between $N$ and  $\mathbb{F}$ as follows:  $\psi^{+} (a)= \alpha $ for any $a\in N$,  $a = \left(\begin{array}{cc}1 & 0\\ \alpha & 1\\ \end{array}\right), \, \, \alpha\in \mathbb{F}$. If $$a_1 = \left(\begin{array}{cc}1 & 0\\ \alpha_1 & 1\\ \end{array}\right)\quad \mathrm{and}\quad a_2 = \left(\begin{array}{cc}1 & 0\\ \alpha_2 & 1\\ \end{array}\right) \, \, \alpha_1, \alpha_2 \in \mathbb{F},$$ then
$$a_1a_2 = \left(\begin{array}{cc}1 & 0\\ \alpha_1 & 1\\ \end{array}\right)\left(\begin{array}{cc}1 & 0\\ \alpha_2 & 1\\ \end{array}\right)=\left(\begin{array}{cc}1 & 0\\ \alpha_1+\alpha_2 & 1\\ \end{array}\right).$$
Thus, $N$ is isomorphic to the additive group $\mathbb{F}^{+}$  which we can consider as a linear space of dimension $r$  over $\mathbb{F}_q$. 

 Define  a bijection $\psi^{\times}$ between $K$ and  $\mathbb{F}\setminus \{0\}$ as follows:  $\psi^{\times} (b)= \beta $ for any $b\in K$, $b = \left(\begin{array}{cc}1 & 0\\ 0 & \beta \\ \end{array}\right), \, \, \beta\in \mathbb{F}\setminus \{0\}$. Clearly, $\psi^{\times}$ is an isomorphism between $K$ and the multiplicative group of $\mathbb{F}$.

Let $K$ be generated by matrix 
$f^{*}  = \left(\begin{array}{cc}1 & 0\\ 0 & \tau \\\end{array}\right)$, where $\tau$ is a primitive element of  $\mathbb{F}$. 

Also let $H$ be a cyclic group generated by $f=(f^{*} )^{q-1} =\left(\begin{array}{cc}1 & 0\\ 0 & \theta \\\end{array}\right)$, where  $\theta=\tau^{q-1}$.  Thus, $G=N H$ is a subgroup of $\mathbb{G} $ of index $q-1$ and the order of $G$ is equal to $q^r (q^r - 1)/(q-1)$. 
Furthermore, $N$ is a normal subgroup of order $q^r$ and index $(q^r - 1)/(q-1)$ in $G$.

By the formula of Gaussian binomial coefficients, the number of $(r-1)$-dimensional subspaces of $N$ equals $t$, where $$t= (q^r - 1)/(q-1).$$ 
Let $\mathbb{M}$ be the set of preimages of all these $(r-1)$-dimensional subspaces of $\mathbb{F}^{+}$ in $N$ under $\psi^{+}$.
Since $\psi^{+}$ is an isomorphism, then the set $\mathbb{M}$ consists
of $t$ subgroups of order $q^{r-1}$  from $N$.

Denote by $M$  one of the subgroups from $\mathbb{M}$.  
Thus, $$\mathbb{M} =\{M_{i}\, \mid\, \psi^{+}(M_i)=\psi^{+}(M)\tau^{i}, \quad i=1,2,\dots ,t\}. $$ 
 
Let $\varphi$ be a permutation on the set $\mathbb{M}=\{M_1,M_2,\dots ,M_t\}$. As usual by $f^i(N\setminus M_{\varphi (i)})$ we denote the set $\{f^i a : a\in N\setminus M_{\varphi (i)}\}.$
\bigskip

{\em Define a subset $S$ of $G$ in the following way:
$$S = \bigcup_{i=1}^{t} f^i(N\setminus M_{\varphi (i)}).$$
}
\bigskip

It is obvious, that $S$ is a generating  set of $G$. In the following lemmas, we examine the question of when this set is closed under inversion.
\medskip

{\bf Lemma 1.}\label{Lemma 1} {\em Subset $S$  is closed under inversion in $G$ if and only if for all integer $i$ the following condition holds}
$$\quad\quad\psi^{+}(M_{\varphi(i)})\theta^{i} = \psi^{+}(M_{\varphi(t-i)}).  \quad\quad\quad\quad(\ast) $$
\medskip

{\bf Proof.}  
Let $s\in S$ and $s=f^i a$, for some integer $i$ and $a\in N\setminus M_{\varphi (i)}$.
Also, let $a=\left(\begin{array}{cc}1 & 0\\ \alpha & 1\\ \end{array}\right)$, for some $\alpha\in \psi^{+}(N\setminus M_{\varphi (i)})$  and $f^i=\left(\begin{array}{cc}1 & 0\\ 0 & \theta^i\\ \end{array}\right)$. 

In such case, $a^{-1}=\left(\begin{array}{cc}1 & 0\\ -\alpha & 1\\ \end{array}\right)$, for  $\alpha\in \psi^{+}(N\setminus M_{\varphi (i)})$  and $f^{-i}=\left(\begin{array}{cc}1 & 0\\ 0 & \theta^{t-i}\\ \end{array}\right)$. 

Then $$f^i a  = \left(\begin{array}{cc}1 & 0\\ 0 & \theta^i\\ \end{array}\right)\left(\begin{array}{cc}1 & 0\\ \alpha & 1\\ \end{array}\right) =
\left(\begin{array}{cc}1 & 0\\ \alpha\theta^i & \theta^i\\ \end{array}\right), $$
$$(f^i a )^{-1} =
\left(\begin{array}{cc}1 & 0\\ -\alpha & 1\\ \end{array}\right) \left(\begin{array}{cc}1 & 0\\ 0 & \theta^{t-i}\\ \end{array}\right)=
\left(\begin{array}{cc}1 & 0\\ (-\alpha\theta^{i})\theta^{t-i} & 1\\ \end{array}\right) \left(\begin{array}{cc}1 & 0\\ 0 & \theta^{t-i}\\ \end{array}\right).$$

Thus $s^{-1}\in S$ if and only if $-\alpha\theta^{i}\in \psi^{+}(M_{\varphi(t-i)})$. 
Hence $S=S^{-1}$ if and only if $$\psi^{+}(M_{\varphi(i)})\theta^{i} = \psi^{+}(M_{\varphi(t-i)}).$$
\medskip

Let $\varphi = (\varphi_1, \dots , \varphi_t)$.
\medskip

{\bf Lemma 2.} {\em  For any $t>2$ there is at least one permutation $\varphi$ satisfying $$\psi^{+}(M_{\varphi(i)})\theta^{i} = \psi^{+}(M_{\varphi(t-i)})$$ 
for all integer $i$. }
\medskip

{\bf Proof.}  Let $\varphi = (\varphi_1, \dots , \varphi_t)$.
 If $t$ is an odd integer, then $$\varphi = (1,t-1,t-3,\dots ,2,t,t-2,t-4,\dots ,3).$$ 

If $t$ is an even integer, then
$$\varphi = (1,t-1,t-3, \dots , t/2+2, t/2-1, t/2-3,\dots , 2,t/2,t-2,t-4,\dots ,$$ $$\dots , t/2+1,t,t/2-2, t/2-4,\dots ,3).$$ 
\bigskip

For example, if $t = 7$, then there are exactly three permutations which satisfy the  condition $(\ast)$ in Lemma 1.
These are $(1, 6, 4, 2, 7, 5, 3)$, $(1, 6, 3, 7, 2, 4, 5)$ and $(1, 6, 2, 3, 5, 7, 4)$.
\bigskip

{\bf Construction 1.} {\em Let $\Gamma$ be a Cayley graph $\mathbb{C}ay(G,S)$ whose vertices are elements of the group $G$ defined as above and  $$S = \bigcup_{i=1}^{t} f^i(N\setminus M_{\varphi (i)})$$.} 
\bigskip

In the next section we prove that if $S$ satisfies the  condition $(\ast)$, then $\Gamma$ is a divisible design graph.  Construction 1 gives us an infinite series of divisible design graphs which are Cayley graphs. Only the first graph among them is known and given in  \cite[Construction 4.20]{HKM}. This divisible design Cayley graph with parameters $(12, 6, 2, 3, 3, 4)$ is the line graph of the octahedron and can be obtained as a Cayley graph from the alternating group of degree $4$ (See Example 1 in Section 4).

\section{Main theorem}

The main goal of our article is to prove the following theorem.
\medskip

{\bf Theorem 2.} {\em Let  $\Gamma$ be a Cayley graph from Construction 1. If $S$ satisfies the  condition $(\ast)$, then $\Gamma$ is a divisible design graph with parameters 
$(v, k, \lambda_1,  \lambda_2, m, n)$, where }
$$v = q^r (q^r - 1)/(q-1),\quad k = q^{r-1}(q^r - 1),$$
$$\lambda_1 = q^{r-1}(q^r - q^{r-1} - 1),\quad \lambda_2 = q^{r-2}(q-1)(q^r - 1),$$ $$m = (q^r - 1)/(q-1),\quad n = q^r.$$
\medskip

{\bf Proof.}
It is clear, that $\Gamma$ is an undirected graph  on $v=q^r (q^r - 1)/(q-1)$ vertices of valency 
$$k=|S| = \Bigl|\bigcup_{i=0}^{t-1}f^i (N\setminus M_{\varphi (i)}) \Bigr| = \sum_{i=0}^{t-1} (q^r - q^{r-1})=$$ $$=(q^r - q^{r-1}) (q^r - 1)/(q-1)= q^{r-1} (q^r - 1).$$ 

Calculate a number of common adjacent vertices for  any two distinct vertices from any coset. Since $\Gamma$ is a Cayley graph, then it is enough to calculate this number for the identity element of $G$ and any non-identity element from $N$.
Let $e$ be the identity element of $G$, $a\neq e$ and $a\in N$. 
It is important to note that $a$ belongs to exactly $t_1 = (q^{r-1} - 1)/(q-1)$ 
 subgroups of $N$ from $\mathbb{M}$. 

Since $\Gamma (e)\cap \Gamma (a)=S \cap S a,$ then 
 $$\Gamma (e)\cap \Gamma (a) =\Bigl( \bigcup_{i=0}^{t-1} f^i(N\setminus M_{\varphi (i)}) \Bigr)\cap \Bigl(\bigcup_{i=0}^{t-1} f^i(N\setminus M_{\varphi (i)}) a \Bigr)=$$ $$= \Bigl(\bigcup_{a\in M_{\varphi (i)}} f^i(N\setminus M_{\varphi (i)})\Bigr) \cup \Bigl(\bigcup_{a\notin M_{\varphi (i)}} f^i(N\setminus M_{\varphi (i)})\cap f^i(N\setminus M_{\varphi (i)}) a \Bigr).$$ Hence $$|\Gamma (e)\cap \Gamma (a)|= \sum_{i=0}^{t_1-1} (q^r - q^{r-1}) + \sum_{i=t_1}^{t-1} (q^r - (q-2)q^{r-1})= q^{r-1}(q^r - q^{r-1} - 1).$$ 

Calculate a number of common adjacent vertices for  any two vertices from any two  distinct cosets. Since 
$\Gamma$ is a Cayley graph, then it is enough to calculate this number for the identity element $e$ from $G$ and any element $g$ from $N f^i$, where $i \neq  0 \pmod{t}$. Let $g=f^x a_g$, $x\neq 0 \pmod{t}$. We have $f^x a_g = [f^{-x}, a_g^{-1}]a_g f^x $, where $[f^{-x}, a_g^{-1}]$ is the commutator of elements $f^{-x}$ and $a_g^{-1}$. It is easy to verify, that if $b_g = [f^{-x}, a_g^{-1}]a_g$, then $b_g \in N$.

Since $ \Gamma (e)\cap \Gamma (g)=S\cap S g$, then 
 $$\Gamma (e)\cap \Gamma (g)=
 \Bigl(\bigcup_{i=0}^{t-1} f^i(N\setminus M_{\varphi (i)})\Bigr) \cap \Bigl(\bigcup_{i=0}^{t-1} f^i(N\setminus M_{\varphi (i)}) g \Bigr).$$

Taking into account that $N$ is a normal subgroup of $G$ and $b_g \in N$ we have  
$$ f^i (N\setminus M_{\varphi (i)})\cap f^j (N\setminus M_{\varphi (j)}) b_g f^x  \neq \emptyset \quad  \mathrm{if\ and\ only\ if} \quad i-x = j. $$   
Thus, we have
 $$ \bigcup_{i=0}^{t-1}\left( f^{i}(N\setminus M_{\varphi (i)}) \cap  
 f^{i-x}( N\setminus M_{\varphi (i-x)}) f^x \right) = 
 \bigcup_{i=0}^{t-1}f^{i-x}\left( f^{x}(N\setminus M_{\varphi (i)}) \cap  ( N\setminus M_{\varphi (i-x)}) f^x \right).$$ 
 
If   $h\in f^{x}(N\setminus M_{\varphi (i)}) \cap  ( N\setminus M_{\varphi (i-x)}) f^x ,$ then there are some  $$\alpha_1 \in \psi^{+}(N\setminus M_{\varphi (i)}), \quad \alpha_2 \in \psi^{+}(N\setminus M_{\varphi (i-x)})$$ such that  
$$h=\left(
\begin{array}{cc}1 & 0\\ 0 &\theta^x\\
\end{array}\right) \left(
\begin{array}{cc}1 & 0\\ \alpha_1 & 1\\
\end{array}\right)= \left(
\begin{array}{cc}1 & 0\\ \alpha_2 & 1\\
\end{array}\right)\left(
\begin{array}{cc}1 & 0\\ 0 & \theta^{x}  \\
\end{array}\right).$$

Therefore, $\theta^x \alpha_1 =  \alpha_2$.

Since bijection $\psi^{+}$ is an isomorphism  between  $N$ and  $\mathbb{F}^+$, then
$$|f^x (N\setminus M_{\varphi (i)}) \cap  
 ( N\setminus M_{\varphi (i-x)}) f^x|= |\psi^{+} (N\setminus M_{\varphi (i)}) \cap  \theta^{-x}
\psi^{+} ( ( N\setminus M_{\varphi (i-x)}))| =$$ 
$$ =|\psi^{+} (N) \setminus \psi^{+} (M_{\varphi (i)}) \cup   \theta^{-x}\psi^{+}( M_{\varphi (i-x)})| = (q^r - 2q^{r-1}+q^{r-2}).$$

Thus, $|\Gamma (e)\cap \Gamma (g)|=$ $$=\sum_{i=0}^{t-1} (q^r - 2q^{r-1}+q^{r-2})=
(q^r - 2q^{r-1}+q^{r-2})(q^r - 1)/(q-1)=$$  $$ = q^{r-2}(q-1)(q^r - 1).$$

Hence, $\Gamma$ is a divisible design graph.  

\section{Examples}
All examples in this section except Example 1 were found using computer search.
\medskip

\noindent{\bf Example 1. Divisible design graph with parameters $(12,6,2,3,3,4)$}.\label{exampl1}

\noindent There is the only example of divisible design Cayley graph  with parameters $(12,6,2,3,3,4)$ basing on the alternating group $\mathrm{Alt}_4$.  We can chose the set $S=\{(13)(24), (12)(34), (123),  (132),  (234), (243)\}$ as its generating set. 
A fragment of the Cayley table of  $\mathrm{Alt}_4$ below shows us the necessary properties.

\noindent\begin{center}
\begin{tabular}[t]{||c|c|c|c|c|c|c|} 
\hline            & (13)(24) & (12)(34) &   (123)  &   (132)  &   (234)  & (243) \\ 
\hline  (13)(24)  &    e     & (14)(23) &   (243)  &   (124)  &   (143)  & (123) \\
\hline  (12)(34)  & (14)(23) &     e    &   (134)  &   (234)  &   (132)  & (142) \\
\hline    (123)   &  (142)   &   (243)  &   (132)  &     e    & (13)(24) & (143) \\
\hline    (132)   &  (234)   &   (143)  &     e    &   (123)  &   (142)  & (12)(34) \\
\hline    (234)   &  (132)   &   (124)  & (12)(34) &   (134)  &   (243)  &    e     \\
\hline    (243)   &  (134)   &   (123)  &   (124)  & (13)(24) &     e    & (234) \\
 \hline \end{tabular}
 \end{center}
\medskip

\noindent{\bf Example 2. Divisible design graphs with parameters $(36,24,15,16,4,9)$}.

\noindent It was found in \cite{GSh} that there are three  non-isomorphic divisible design graphs with parameters $(36,24,15,16,4,9)$. From our Construction 1 with permutations $(1,3,4,2)$ and $(1,3,2,4)$ we have two of that non-isomorphic divisible design graphs, which are based on subgroup of index 2 of $AG(9)$. It is important to note that, if $t$ is even, then  $t/2$ and $t$ can be in any place in $varphi$  according to the condition $(\ast)$.

\medskip

\noindent{\bf Example 3. Divisible design graphs with parameters $(56,28,12,14,7,8)$}.

\noindent It was found in \cite{GSh} that there are five non-isomorphic examples of divisible design graphs  with parameters 
$(56,28,12,14,7,8)$
 which are based on group $AG(8)$. This group can be described as follows $$G=\left<\,f_1,\, f_2,\, f_3,\, f_4\,\right> $$ with defining group relationships $$ f_1^7= f_2^2= f_3^2= f_4^2=e,$$ $$\, f_2*f_1=
f_1*f_3,\, f_3*f_1=f_1*f_4,\,
 f_4*f_1= f_1*f_2*f_4.$$ 
Below, we give the list of generating sets for these Cayley graphs. 
\medskip

\noindent$S_1 =\{f_3, f_2*f_3, f_3*f_4, f_2*f_3*f_4,  f_1*f_4, f_1*f_2*f_4, f_1*f_2*f_3*f_4, f_1*f_3*f_4, f_1^2*f_2, f_1^2*f_3, f_1^2*f_4, f_1^2*f_2*f_3*f_4, f_1^3*f_2, f_1^3*f_3, f_1^3*f_2*f_4, f_1^3*f_3*f_4, f_1^4*f_2, f_1^4*f_2*f_3, f_1^4*f_2*f_4, f_1^4*f_2*f_3*f_4,  f_1^5*f_2, f_1^5*f_4,  f_1^5*f_2*f_3, f_1^5*f_3*f_4, f_1^6*f_3, f_1^6*f_4, f_1^6*f_2*f_3, f_1^6*f_2*f_4\};$
\medskip

\noindent$S_2 =\{ f_2, f_4, f_2*f_3, f_3*f_4, 
f_1, f_1*f_3, f_1*f_4, f_1*f_3*f_4, 
f_1^2*f_3, f_1^2*f_4, f_1^2*f_2*f_3,  f_1^2*f_2*f_4,  f_1^3, f_1^3*f_2*f_4, f_1^3*f_3*f_4, f_1^3*f_2*f_3,  f_1^4,  f_1^4*f_2, f_1^4*f_4, f_1^4*f2*f_4, f_1^5*f_2, f_1^5*f_3, f_1^5*f_2*f_4, f_1^5*f_3*f_4, f_1^6,   f_1^6*f_2, f_1^6*f_3,  f_1^6*f_2*f_3\};$
\medskip

\noindent$S_3=\{f_1, f_3, f_4, f_1^2, f_1*f3, f_2*f_4, f_3*f4, f_1^2*f_2, f_1^2*f_3, f_1^2*f_4, f_1*f_2*f_4, f_1*f_3*f_4, f_1^3*f_3,f_1^5, f_1^4*f_2, f_1^4*f_4, f_1^3*f_2*f_3, f_1^3*f_2*f_4, f_1^3*f_3*f_4, f_1^6, f_1^5*f_2, f_1^4*f_2*f_3,f_1^4*f_2*f_4, f_1^6*f_2, f_1^6*f_4, f_1^5*f_2*f_3, f_1^5*f_3*f_4, f_1^6*f_2*f_3\};$
\medskip

\noindent$S_4=\{f_1, f_2, f_3, f_1*f_3, f_1*f_4, f_2*f_4, f_3*f_4, f_1^2*f_2, f_1^2*f_3, f_1^2*f_4, f_1*f_3*f_4, f_1^3*f_3, f_1^4*f_3,f_1^4*f_4, f_1^3*f_2*f_3, f_1^3*f_3*f_4, f_1^2*f_2*f_3*f_4, f_1^6, f_1^5*f_2, f_1^5*f_4, f_1^4*f_2*f_3,f_1^4*f_2*f_4, f_1^3*f_2*f_3*f_4, f_1^6*f_2, f_1^6*f_3, f_1^5*f_2*f_3, f_1^5*f_3*f_4, f_1^6*f_2*f_3\};$
\medskip

\noindent$S_5=\{f_1, f_2, f_3, f_4, f_1^2, f_1*f_3, f_1*f_4, f_1^2*f_4, f_1*f_3*f_4, f_2*f_3*f_4, f_1^3*f_3, f_1^3*f_4, f_1^2*f_2*f_3,f_1^5, f_1^4*f_2, f_1^4*f_4, f_1^3*f_2*f_3, f_1^3*f_2*f_4, f_1^2*f_2*f_3*f_4, f_1^6, f_1^5*f_2, f_1^5*f_4,f_1^4*f_2*f_3, f_1^4*f_3*f_4, f_1^6*f_2, f_1^6*f_3, f_1^5*f_2*f_4, f_1^6*f_2*f_3\}.$
\medskip

Three of them are isomorphic to the graphs we have from our Construction 1 with permutations which pointed out after Lemma 2.
\medskip

\noindent{\bf Example 4. Divisible design graph with parameters $(80,60,44,45,5,16)$}.

\noindent  There is at least one divisible design Cayley graph which is based on subgroup of index 3 of $AG(4^2)$ that we have from our Construction 1 with permutations $(1,4,2,5,3)$. This is the first example where $q$ is not prime.

\section{Conclusion remarks}
Any divisible design graph can be considered as a symmetric group-divisible design, the vertices of which are points, and the neighborhoods of the vertices are blocks. Such a design is called the neighborhood design of a graph.  However, non-isomorphic graphs can correspond to isomorphic designs. Examples $2$ and $3$ of this article give us non-isomorphic divisible design graphs which produce isomorphic group divisible designs. If group-divisible design admit a symmetric incidence matrix with zero diagonal, then it corresponds to divisible design graph  (see Section 4.3. in \cite{HKM}).

There is a great possibility to construct divisible designs from groups.
 Let $G$ be a group of order $mn$ containing a subgroup $N$ of order $n$. A $k$-subset $D$  of $G$  is called a {\em divisible} $(m, n, k, \lambda_1,  \lambda_2)$ difference set ({\em divisible by cosets of subgroup } $N$ ) if the list of elements $xy^{-1}$ with $x, y \in D$ contains all non-identity elements in $N$ exactly $\lambda_1$ times and all elements in $G \setminus N$ exactly $\lambda_2$ times. In case that $N = \{0\}$, the definition of a divisible difference set coincides with the definition of a {\em difference set} in the usual sense. In case that $\lambda_1=0$, the definition of a divisible difference set coincides with the definition of a {\em relative difference set} \cite{BJL, P1995, P1996}.

Divisible difference set $D$ gives rise to a symmetric group-divisible design
$\mathbb{D}$ with the set of blocks $\{D g |\, g\in G\}$ and has the same parameters as $D$. This symmetric group-divisible design is called the {\em  development of $D$} and admits $G$ as a regular automorphism group (by right translation). Thus, symmetric group-divisible designs with a regular group $G$ are equivalent to divisible difference sets in $G$. For having
a symmetric incidence matrix with zero diagonal, the divisible difference set should be reversible
(or equivalently, it must have a strong multiplier $-1$). There is more information on such difference sets in \cite{AJP}.

\section*{Acknowledgements} The authors are grateful to Vladimir Trofimov, Sergey Goryainov, Elena Konstantinova and participants of workshop "New trends in algebraic graph theory", which was organized by  Mathematical Center in Akademgorodok, for useful discussions.

\bigskip

\noindent{\bf ORCID}
\medskip

\noindent Vladislav V. Kabanov http://orcid.org/0000-0001-7520-3302

\noindent Leonid Shalaginov  https://orcid.org/0000-0001-6912-2493
\bigskip


\end{document}